\documentstyle{article}
\begin{document}
\title{{\bf
 Nonclassical Solutions of Fully Nonlinear Elliptic
Equations  } }

\author{{Nikolai Nadirashvili\thanks{LATP, CMI, 39, rue F. Joliot-Curie, 13453
Marseille  FRANCE, nicolas@cmi.univ-mrs.fr},\hskip .4 cm Serge
Vl\u adu\c t\thanks{IML, Luminy, case 907, 13288 Marseille Cedex
FRANCE, vladut@iml.univ-mrs.fr} }}

\date{}
\maketitle
\def\n{\hfill\break}
\def\al{\alpha}
\def\be{\beta}
\def\ga{\gamma}
\def\Ga{\Gamma}
\def\om{\omega}
\def\Om{\Omega}
\def\ka{\kappa}
\def\lm{\lambda}
\def\Lm{\Lambda}
\def\dl{\delta}
\def\Dl{\Delta}
\def\vph{\varphi}
\def\vep{\varepsilon}
\def\th{\theta}
\def\Th{\Theta}
\def\vth{\vartheta}
\def\sg{\sigma}
\def\Sg{\Sigma}

\def\bendproof{$\hfill \blacksquare$}
\def\wendproof{$\hfill \square$}

\def\holim{\mathop{\rm holim}}
\def\span{{\rm span}}
\def\mod{{\rm mod}}
\def\rank{{\rm rank}}
\def\bsl{{\backslash}}

\def\il{\int\limits}
\def\pt{{\partial}}
\def\lra{{\longrightarrow}}

{\em Abstract.} We prove the existence of non-smooth solutions to
fully nonlinear elliptic equations.

\section{ Introduction}
\bigskip

In this paper we study the regularity of solutions of fully
nonlinear elliptic equations of the form
$$F(D^2u)=0\leqno(1)$$
defined in a domain of ${\bf R}^n$. Here $D^2u$ denotes the
Hessian of the function $u$. We assume that
   $F$ is uniformly elliptic, i.e. there exists a constant $\Lm\ge 1$
(called an {\it ellipticity constant\/}) such that
$$\Lm^{-1}|\xi|^2\le
F_{u_{ij}}\xi_i\xi_j\le\Lm |\xi |^2\;,\forall\xi\in {\bf
R}^n\;.\leqno(2)$$ Here, $u_{ij}$ denotes the partial derivative
$\pt^2 u/\pt x_i\pt x_j$. A function $u$ is called a {\it
classical\/} solution of (1) if $u\in C^2(\Om)$ and $u$ satisfies
(1).  Actually, any classical solution of (1) is a smooth
($C^{\alpha +3}$) solution, provided that $F$ is a smooth
$(C^\alpha )$ function of its arguments.

Consider the following Dirichlet problem
$$\cases{F(D^2u)=0 &in $\Om$\cr
u=\vph &on $\pt\Om\;,$\cr}\leqno(3)$$ where  $\Omega \subset {\bf
R}^n$ is a bounded domain with smooth boundary $\partial \Omega$
and $\vph$ is a continuous function on $\pt\Om$.

It is not difficult to prove that problem (3) has no more than one
classical solution (see e.g. [GT]). The basic problem is  the
existence of such classical solutions. Although the first
systematic  study of the Dirichlet problem for fully nonlinear
equations was done by Bernstein at the beginning of the 20-th
century (see [GT]), the first complete result didn't appear until
1953, when Nirenberg proved the existence of a classical solution
to problem (3) in dimension $n=2$ ([N]). For $n\ge 3$, the problem
of the existence of classical solutions to Dirichlet problem (3)
remained open.

In order to get a solution to the problem (3) one can try to
extend the notion of the classical solution of the equation (1).
That was done recently: Crandall-Lions and Evans developed the
concept of viscosity (weak) solutions of the fully nonlinear
elliptic equations. As a characteristic property for such
extension can be taken the maximum principle in the following
form:

Let $u_1,u_2$ be two  solutions of the following equations,
$F(D^2u_1)=f_1$ in $\Omega$ and $F(D^2u_2)=f_2$ in $\Omega$. Then
for any subdomain $G\subset \Omega$ the inequalities $f_1\leq f_2\
(f_1\geq f_2)$ in $G$ and $u_1\geq u_2\ (u_1\leq u_2)$ on
$\partial G$ imply the inequality $u_1\geq u_2\ (u_1\leq u_2)$ in
$G$.

Such maximum principle holds for $C^2$ functions $u_1,u_2$. We
call a continuous function $u_1$ a viscosity solution of
$F(D^2u_1)=f_1$ if the above maximum principle holds for $u_1$ and
all $C^2$-functions $u_2$.

It is possible to prove the existence of a viscosity solution to
the Dirichlet problem (3) and Jensen's theorem says that the
viscosity solution of the  problem (3) is unique. For more details
see [CC], [CIL].

There are important classes of the fully nonlinear Dirichlet
problems for which the viscosity solution is in fact a classical
one, e.g., due to Krylov-Evans regularity theory, in the case when
the function $F$ is convex, (see [CC], [K] ). However, for the
general $F$ the problem of the coincidence of viscosity solutions
with the classical remained open.

The central result of this paper is the existence of nonclassical
viscosity solution of (1) in the dimension 12. More precisely we
prove

\medskip

{\bf Theorem.} {\it The function
$$w(x):={Re (\omega_1\omega_2\omega_3)\over |x|}, $$
where $\omega_i\in {\bf H},\ i=1,2,3,$ are Hamiltonian
quaternions, $x=(\omega_1,\omega_2,\omega_3)\in {\bf H}^3={\bf
R}^{12}$ is a viscosity solution in ${\bf R}^{12}$ of a uniformly
elliptic equation $(1)$ with a smooth $F$.}
\medskip

One can find the explicit expression for $w$ in the coordinates of
${\bf R}^{12}$  in Sections 3 and 4. The elliptic operator $F$
will be defined in a constructive way in Section 2, and its
ellipticity constant $\Lm <  10^{8}$.

As an immediate consequence of the theorem we have

\medskip

{\bf Corollary.} {\it Let $\Om \subset {\bf R}^{12}$ be the unit
ball and $\varphi = w$ on $\pt \Om$. Then there exists a smooth
uniformly elliptic  $F$ such that the Dirichlet problem $(2)$ has
no classical solution.}
\medskip

Homogeneous order 2 function $w$ is smooth in ${\bf
R}^{12}\setminus \{0\}$ and has discontinuous second derivatives
at $0$. It is interesting to notice that the set of homogeneous
order $\alpha\in{\bf R} $ solutions of (1) in ${\bf R}^n\setminus
\{0\}$ for $\alpha\neq 2$ has a simple structure: each such
solution of (1) has to be also a solution of a linear elliptic
equation with constant coefficients, [NY].

The question on the minimal dimension $n$ for which there exist
nontrivial homogeneous order 2 solutions of (1) remains open. We
notice that from the result of Alexsandrov [A] it follows that any
homogeneous order 2 solution of the equation (1) in ${\bf R}^3$
with a real analytic $F$ should be a quadratic polynomial. For a
smooth and less regular $F$ similar results in the dimension 3 one
can find in [HNY].

\medskip
{\em Acknowledgements.} The authors would like to thank S. Kuksin
and Y.Yuan for very useful discussions as well as the anonymous
referee for his very pertinent remarks.

\section{ The
Hessian Problem}

\subsection{}

Let $w$ be a homogeneous function of order 2, defined on ${\bf
R}^n$ and smooth in ${\bf R}^n\setminus\{0\}$. Then the Hessian of
$w$ is homogeneous of order 0, and defines a map
$$H:S^{n-1} \to D^2w \in Q\;,$$
where $Q$ denotes the space of quadratic forms on ${\bf R}^n$,
which we will sometimes identify as symmetric $n\times n$
matrices, $Q\simeq S^{n\times n}$. The inner product of $a,b \in
Q$ is given by $a\cdot b = trace (ab)$. We say that $w$ {\it
satisfies property \/} (H) ($w$ is a solution of the Hessian
Problem) if the following holds:

\medskip\noindent
{(H)} 1).{\it The map $H:S^{n-1} \to  Q$ is a smooth embedding.}

2). {\it  There exists a constant $M\ge 1$ such that for any two
points $a,b \in H(S^{n-1}), a\not= b$, if $\mu_1 \ge \cdots \ge
\mu_n$ denote the eigenvalues of the quadratic form $a-b$, then}
$${1/ M} <-{\mu_1/\mu_n}< M\;.\leqno
(4)$$
\smallskip\noindent

\medskip\noindent
{\bf Main Lemma.} {\it If function $w$ satisfies hypotheses} (H)
{\it then $w$ is a viscosity solution in ${\bf R}^n$ of a
uniformly elliptic equation $(1)$ .}
\bigskip\noindent
\subsection{}

 Let us choose in the
space

$Q$ an orthonormal coordinate system $z_1,\dots,z_k,s,$
$k={n(n+1)\over 2} -1$ such that $\sqrt ns$ is the trace. Let $\pi
: Q\to Z$ be the orthogonal projection of $Q$ onto the $z$-space.
For $\lambda \ge 1$, we denote by $K_\lambda$ the cone
 $$
 K_\lambda=
\{a\in Q: {\hbox{ there exists}}\; C>0 {\hbox{ s. t. the
eigenvalues of }} a \in [C/\lambda, C\lambda] \;\}.
$$
Notice, that inequalities (2) mean that the eigenvalues of $\nabla
F$ are on the segment $[\Lambda^{-1},\Lambda ]$. In particular (2)
implies the inclusion $\nabla F\in K_\Lambda$. Since $0\notin
K_\Lambda$ it follows in particular that $F^{(-1)}(0)\in
C^\infty$.

Since on $Q$ the maximal eigenvalue of a quadratic form is a
convex function and the minimal eigenvalue is a concave function
it follows that $K_\Lambda$ is a convex cone.

Let $K_\lambda^\ast$ denote the adjoint cone of $K_\lambda$, that
is,
$$
K_\lambda^\ast = \{b\in Q: b\cdot c \geq 0 {\hbox{ for all }}\; c
\in K_\lambda \}\; .
$$

As an adjoint to a convex cone the cone $K_\lambda^\ast$ is a
convex cone itself.

Set $L_\lambda \: =Q \setminus
(K_\lambda^\ast\cup-K_\lambda^\ast)\;.$ Notice that $a \in
L_\lambda $ is equivalent to $a\cdot b =0$ for some $b \in
K_\lambda$, i.e., $L_\lambda$ is a union of all hyper-planes in
$Q$ with normals in $K_\lambda$. Since the quadratic forms of
$K_\lambda$ are positively defined it follows that the vector
$I\in K_\lambda^\ast$. Let $K\subset Q $ be a cone with a smooth
strictly convex base such that $K_{2\lambda}\subset K\subset
K_\lambda .$
  Let $e_1,...e_k,\ I$ be an
orthonormal basis of $Q$ corresponding to the coordinates
$z_1,...,z_k,\ s$. Then any matrix $b\in Q$ can be written as
$$b\ =\ sI\ +\ \sum_{j=1}^k z_je_j.$$
Now define
$$x(z)\ :=\inf\{ c:\ (a+cI)\in K^* \}$$
for $a:= \sum_{j=1}^k z_je_j$.  The graph of the function $s=x(z)$
represents the boundary of the cone $K$. Clearly $x(\cdot )$ is
Lipshitz, convex, homogenous, smooth outside the origin and
$x(0)=0$.  By a simple computation we get that $|\nabla x|<\sqrt
n$.

\smallskip
Let $G\subset Q$ be a set. We say that $G$ {\it satisfies the
$\lambda$-cone condition} if for any two points $a,b \in G$, the
matrix $a-b \in L_\lambda.$
\smallskip

{\bf Lemma.} {\it Let $\Sigma \subset Q$ be a smooth compact
$(n-1)$-dimensional manifold. Assume that $\Sigma $ satisfies the
$\lambda$-cone condition. Then there exists a smooth function $F$
on $Q$ such that $F(\Sigma)=0$, and which satisfies the inequality
(2) with the ellipticity constant $\Lambda <4\lambda^2\sqrt n$.}

{\it Proof of the lemma.} Set $\sigma=\pi \bigl (\Sigma )$. We
prove that $\Sigma$ is a graph of a Lipschitz continuous function,
$$\Sigma =\{z\in
\sigma:s=g(z)\}\;.$$ Let $a,\hat a \in \Sigma , \ a =\ sI+
\sum_{j=1}^k z_je_j, \ \hat a\ =\ \hat sI+\sum_{j=1}^k\hat
z_je_j$. Since $a-\hat a \in L_\lambda$, we have $-x(z-\hat z) \le
\hat s - s \le x(z- \hat z)$. Since $x(0)=0,\ \ g(z):=s$ is
single-valued. Also \smallskip $|g(z)-g(\hat z)|=|s-\hat s|\le
|x(z-\hat z)|\le C|z-\hat z|.$

Hence, since $\Sigma$ is a smooth surface $g$ is a smooth function
and $\sigma $ is a smooth surface as well.

Let $G^m_k$ be the Grassmannian manifold of $m$-dimensional
subspaces of the $k$-dimensional subspace $z$ of $Q$. Let $l\in
G^{n-1}_k$ and $t:l\rightarrow s$ be a linear function on $l$,
such that the graph of $t$ satisfies the $\lambda$-cone condition.
All such linear functions $t$ defined on all $l\in G^{n-1}_k$ we
denote by $\tau$.  Let $t\in \tau$ defined on $l\in G^{n-1}_k$ be
such that $|\nabla t|\neq0$. Then there exists a constant $c'>1$
such that $c't\leq x$ on $l$ and there is a point $z'\in l,\
|z'|=1$ with $c't(z')=x(z')$.  Since $K^\ast$ is a strictly convex
cone the vector $z'$ is unique. Denote
$$\eta(t)=\{z\in l , t(z)=0\}.$$
Then $\eta (t) \perp \nabla t$. Since $\eta (t) \subset Z $ and
$\eta (t)$ is tangent to the cone $K^*$ at $z', x(z')$ it follows
that $\eta (t) \bot \nabla x(z')$.

Let $\theta $ be a smooth function defined on $[1, \infty)$ such
that $0\leq \theta \leq 1, \ \theta =1$ on $[1,A]$, $\theta =0$ on
$[2A,\infty)$, where a sufficiently large constant $A$ will be
chosen later. Set
$$\nu (t)= \theta(c')\nabla x (z')+ (1-\theta(c'))\nabla
t.$$

For $z\in \sigma$ we denote by $l(z)\in G^{(n-1)}_k$ the tangent
subspace to $\sigma $ at $z$. Let $t_z(x),\ x\in l(z)$, be the
differential of $g$ at $z$.

Let $z\in \sigma$. Denote by $\Psi(z)$ $(n-1)$-dimensional
subspace spanned by $\eta (t_z)$ and $\nu (t_z)$, if $\nabla
t_z\neq 0$. If $\nabla t_z= 0$ we set $\Psi(z)=l$. Thus we defined
a smooth map
$$\Psi :\sigma \rightarrow
G^{n-1}_k.$$

There exists $\gamma \subset Z$ a closed neighborhood of $\sigma$
such that $\gamma $ is diffeomorphic to $\sigma \times B$, where
$B$ is the $(k-n+1)$-dimensional disk. We define a  projection
$\gamma \to \sigma $ such that the fiber $p^{-1}(z)$ is orthogonal
to $\Psi (z)$ at $z\in \sigma$. Since $(\nabla x(z'), z')>0$ the
fiber $p^{-1}(z)$ is transversal to $\sigma $ at $z$.
 We extend the function $g$  to
$\gamma $ by $g(y)=g(p(y))$. Let $\Gamma $ be the graph of $g$
over $\gamma $. Let $z\in \sigma $ and $dg(z)$ be the differential
of $g$ over $\gamma $. For sufficiently large constant $A$ the
following alternative holds: either $|\nabla g(z)|$ is
sufficiently small, or the graph of $c'dg(z)$ is tangent to the
cone $K^*$. In  both cases the graph of $dg(z)$ satisfies
$2\lambda$-cone condition. Since $g\in C^1(Z)$ and along the
fibers the function $g$ is a constant, we may assume the
neighborhood $\gamma $ to be sufficiently small so that $\Gamma $
satisfies the $\lambda $-cone condition.

Since $K_{2\lambda}\subset K\ g\in C^1(\sigma )$ the function $g$
has an extension $\widetilde{g}$ from the set $\gamma$ to ${\bf
R}^k$ such that $\widetilde{g}$ is a Lipschitz function and the
graph of $\widetilde{g}$ satisfies the $2\lambda$-cone condition.
One can define such an extension $\widetilde{g}$ simply by the
formula
$$\widetilde{g}(z):=\inf_{w\in \gamma } \bigl\{ g(w)\ +\
x(z-w) \bigr\}\;.$$

To demonstrate that this formula works let $(z, \tilde g(z)), \
(\hat z,\tilde g(\hat z))$ lie on the graph $\tilde g$. We must
show
$$-x(z-\hat z)\ \le \ \tilde g(z)-\tilde g(\hat z)\ \le x(z-\hat
z).$$ Now
$$\tilde g(\hat z)\ =\ g(w)+x(\hat z-w)$$
for some $w\in \gamma$. Thus
$$\tilde g(z)-\tilde g(\hat z)\ \le \ g(w(\hat
z))+x(z-w(\hat z))-(g(w(\hat z))+x(\hat z-w(\hat z)))\ \le \
x(z-\hat z),$$ since $x(a+b)\ \le \ x(a)+x(b)$, as $x(\cdot )$ is
convex and homogenous. Similarly
$$\tilde g(z)-\tilde g(\hat z) \ \ge \ -x(z-\hat
z).$$
\medskip
Let $D_1 , D_2 \subset Z$ be bounded domains such that $\sigma
\subset D_1 \subset  D_2 \subset \subset \gamma.$ Next, let $l(z)
\in C^\infty({{\bf R}^k})$ be supported on the unit ball
$B_1\subset {\bf R}^k$ and $\int_{{\bf R}^k}l(z)\,dz=1$ and set
$$l_\delta(z)={1\over\delta^k}\,l\left({z\over\delta}\right).$$
Let $h \in C^{\infty}(Z)$,$h=1$ on $D_1$, $h=0$ on $Z\setminus
D_2, 0\leq h \leq 1$ on $Z$.
\medskip
Set
$$g_\epsilon = \widetilde g \ast l_\epsilon
,$$
$$G_\epsilon = hg + (1-h)g_\epsilon.$$

Since the graph of the function $g$ satisfies the $\lambda $-cone
condition it follows that the upper normals to the graph is in the
cone $K$. Since $K^\ast$ is a convex cone the upper normals to the
graphs of the functions $g_\epsilon$ satisfies the $2\lambda$-cone
condition for all small $\epsilon >0$, and hence the graphs of
linear function $ d_z g_\epsilon$ is in $L_{2\lambda} $ for all
$z$ where $d_z$ is the differential at $z$. Since the functions
$g_\epsilon$ are defined on the whole space ${\bf R}^k$ it follows
that the graphs of the functions $g_\epsilon$ satisfies the
$\lambda$-cone condition. Really, let $ a,b,\; a\neq b$  be on the
graph of $g_\epsilon$. If $ a-b\notin L_{2\lambda}$ then there is
a point  $\alpha \in [\pi(a), \pi(b)]$ such that $ d_\alpha
g_\epsilon \notin L_{2\lambda}.$ For any $k>0$ the function $
g_\epsilon \to g$ in $C^k( D_2)$ as $\epsilon \to 0$. Hence for a
sufficiently small $\epsilon_o >0$ the graph of the function
$G_{{\epsilon}_0} :=G$ will satisfy the $2\lambda$-cone condition.
 Moreover $G$ will be a smooth function on $Z$,
$G=g$ on $D_1$ and $|\nabla G|<\sqrt n$ on $Z$.
\medskip

Let us set
$$
F:= s - G(z).
$$
Denote
$$b\ :=\ \nabla   F  \ = \ (-\nabla G, 1)\ ,$$
$$a\ :=\
(\nabla  G / |\nabla  G|,\ |\nabla  G| ).$$ The vector $a$ is
tangent to the level surface of the function $ F$, and  $tr( b)=
\sqrt n$. Since level surfaces of the function  $  F$ satisfies
the $2\lambda$-cone condition and $a\cdot b = 0$, it follows that
$a \in L_{2\lambda}$ and hence $b\in K_{2\lambda}$.  Therefore the
function $F$ satisfies the ellipticity conditions with the
ellipticity constant $\Lambda < 4\lambda^2\sqrt n$.
\medskip

{\bf Remark 1.} For a real-analytic manifold $\Sigma $ one can
obtain the existence of a real-analytic function $F$ after
insignificant changes in the construction.
\medskip

{\bf Remark 2.} The proof of the lemma holds if instead of
compactness of $\Sigma $ we assume that $\Sigma$ is a smooth
closed manifold with a boundary.
\medskip

{\bf Proof of Main Lemma.} Set
$$\lambda = (n-1)M$$
 Let $\xi $ and $\eta $ be correspondingly
negative and nonnegative subspace of the quadratic form $a-b$ in
${\bf R}^n$. Denote by $c\in Q$ the quadratic form $l|\xi |^2 +
m|\eta |^2,\ l,m>0$ such that $(a-b)\cdot c = 0$. Then
$1/(n-1)M<l/m<(n-1)M$ by (H) 2) and hence  the set $H(S^{n-1})$
satisfies the $\lambda $-cone condition.

For $\Sigma = H(S^{n-1})$ we define function $F$ by Lemma. Then
the function $w$ satisfies the equation
$$
  F(D^2 w)= 0.
$$
on ${\bf R}^n\setminus \{0\}$.
\medskip
We show now that $w$ is a viscosity solution of (1) on the whole
space ${\bf R}^n$.

Let $p(x),\ x\in {\bf R}^n$ be a quadratic form such that $p\leq
w$ on ${\bf R}^n$. We choose any quadratic form $p'(x)$ such that
$p\leq p'\leq w$ and there is a point $x'\neq 0$ at which $p'(x')=
w(x')$. Then it follows that $F(p)\leq F(p')\leq 0$. Consequently
for any quadratic form $p(x)$ from the inequality $p\leq w$
($p\geq w$) it follows that $F(p)\leq 0$ ($F(p)\geq 0$). This
implies that $w$ is a viscosity solution of (1) in ${\bf R}^n$
(see Proposition 2.4 in [CC]).

\section{Cubic form $P $ }In this section we introduce and
investigate the cubic form which will be used to construct our
non-classical solutions. Let $ V=(X,Y,Z)\in {\bf R}^{12}$ be a
variable vector with $X,Y,$ and $Z\in {\bf R}^4.$ For any $
t=(t_0,t_1,t_2,t_3)\in {\bf R}^4$ we denote  by $ qt=t_0+t_1\cdot
i+t_2\cdot j+t_3\cdot k\in {\bf H} $ (Hamilton quaternions). For
any $ q=q_0+q_1\cdot i+q_2\cdot j+q_3\cdot k\in {\bf H} $ its
conjugate will be denoted $ q^*=q_0-q_1\cdot i-q_2\cdot j-q_3\cdot
k$; thus, $q^*q=q q^*=\mid q\mid^2=q_0^2+q_1^2+q_2^2+q_3^2.$

    Define the cubic form $P=P(V)=P(X,Y,Z) $ as follows
    $$ P(X,Y,Z)=Re(qX\cdot qY\cdot qZ)=X_0Y_0Z_0-
    X_0Y_1Z_1-X_0Y_2Z_2-X_0Y_3Z_3$$
    $$-X_1Y_0Z_1- X_1Y_1Z_0-X_1Y_2Z_3+X_1Y_3Z_2-
    X_2Y_0Z_2+X_2Y_1Z_3-X_2Y_2Z_0-X_2Y_3Z_1 $$
    $$ -X_3Y_0Z_3- X_3Y_1Z_2+X_3Y_2Z_1-X_3Y_3Z_0.$$

    Let $d=(a,b,c)\in {\bf R}^{12} $ be a vector with the norm
    $ \sqrt 3$,

    $\mid\mid a\mid\mid^2+\mid\mid b\mid\mid^2+\mid\mid c\mid\mid^2=3.$
    Define the quadratic form $$Q_d=Q_{a,b,c}= Q_{a,b,c}(X,Y,Z)$$
    by differentiating $ P$ in the direction $ d$:

    $$
    Q_{a,b,c}(X,Y,Z)=\sum_{i=0}^4 a_i \pt  P/\pt X_i+\sum_{i=0}^4 b_i
\pt  P/\pt Y_i+
    \sum_{i=0}^4 c_i \pt  P/\pt Z_i\; .
    $$

    A direct calculation shows that
    $$ Q_d(X,Y,Z)=X^t M_c Y+X^t M_b^t Z+Y^t M_a Z$$
    where, in general, we define the matrix $M_s $ for an arbitrary
$s\in {\bf R}^4$ by $$M_s=
  \left(%
\begin{array}{cccc}
  s_0&-s_1&-s_2&-s_3 \\
   -s_1&-s_0&-s_3&\; s_2 \\
  -s_2& s_3&-s_0&-s_1 \\
   -s_3&-s_2& s_1&-s_0\\
\end{array}%
\right)$$ Direct (and easy) calculations show that $M_s $ has the
following properties:

\medskip
1).$\quad\quad M_s\cdot M_s^t=M_s^t \cdot M_s= \mid\mid
s\mid\mid^2 I_4;$

\noindent thus, $M_s $  is proportional to an orthogonal matrix.
In particular, if  $\mid\mid s\mid\mid=1$    then $M_s $  is
    orthogonal itself. In general, we write $M_s= \mid\mid
s\mid\mid O_s$  with $O_s\in O(4).$

\medskip
2).  $\det (M_s)=-\mid\mid s\mid\mid^4,$  $\quad\quad \det
(O_s)=-1;$

\medskip
3). the characteristic polynomial $ PM_s(x)$ of $M_s $  factors as
$$ PM_s(x)=(x^2-\mid\mid s\mid\mid^2)(x^2+2 s_0 x +\mid\mid
s\mid\mid^2)$$ with $s_0=Re(qs); $ and that of $ O_s$   as
$$
PO_s(x)=(x^2-1)(x^2+2 s^*_0 x +1)$$  with $s^*_0=s_0/\mid
qs\mid=Re(qs/\mid qs\mid); $

\medskip
4). define the symmetric matrix $N_s=  (O_s+ O_s^t)  ;$  then its
characteristic polynomial $ PN_s(x)=(x^2-4)(x+2 s^*_0    )^2$, its
spectrum  being
     \\$ Sp(N_s)=\{ 2,-2,-2 s^*_0,-2 s^*_0\}$;\\

     \medskip
5). $M_s $ is the matrix (with respect to the standard basis) of
the endomorphism ${\bf H}\rightarrow {\bf H} , $ $ q\mapsto
\bar{q}{\cdot\bar qs}.$
 \medskip

  The points 3 and 5 applied to the product matrix  $
M_{rst}=M_r\cdot M_s\cdot M_t$, $r,s,t $ being arbitrary vectors
in  ${\bf R}^4$ give the  following formula for the characteristic
polynomial  $PM_{rst}$  of $ M_{rst}$:
$$\quad PM_{rst} (x)  =(x^2-\mid\mid r\mid\mid^2
  \mid\mid s\mid\mid^2\mid\mid t\mid\mid^2)
(x^2+2P(r,s,t)x+\mid\mid r\mid\mid^2 \mid\mid s\mid\mid^2\mid\mid
t\mid\mid^2)
$$
with $P(r,s,t)=Re(qr \cdot qs\cdot qt)$ as above. Indeed, $
M_{rst}$ is conjugate to the matrix of the endomorphism $ q\mapsto
\bar{q}{\cdot\bar qr}{\cdot\bar qs}{\cdot\bar qt}.$

 \par       For
the corresponding orthogonal matrix $ O_{rst}$  we get the
polynomial
$$   PO_{rst} (x)=(x^2-1)(x^2+2\bar P(r,s,t)x+1)$$ where
   $\bar P(r,s,t)=P(r,s,t)/(\mid\mid r\mid\mid\cdot
\mid\mid s\mid\mid \cdot \mid\mid t\mid\mid)$      and for the
corresponding symmetric  matrix $N_{rst}=O_{rst}+O_{rst}^t$ the
spectrum   is
  $$(*)\quad Sp(N_{rst} )=\{ 2,-2, \;-2\bar P(r,s,t),-2\bar P(r,s,t)\}.
  $$\par

{\em Warning}: in the case of the product  of two matrices $M_r,
M_s$ the characteristic  polynomial is completely different;
namely, if $ M_{rs }=M_r\cdot M_s $ then $$ PM_{rs} (x)=
   (x^2-2(r,s)x+\mid\mid r\mid\mid^2 \mid\mid
s\mid\mid^2)^2$$ with the usual scalar product $(r,s)$. \par

Define now two quantities $ m= m(d)=m(a,b,c)=\mid\mid
a\mid\mid\cdot \mid\mid b\mid\mid \cdot \mid\mid c\mid\mid $, $
n=n(d)=n(a,b,c)=P(a,b,c)$. Clearly, $\mid n(d)\mid\le m(d) \le 1$
by  the  inequality between the geometric and quadratic means,
since $\mid\mid a\mid\mid^2+\mid\mid b\mid\mid^2+\mid\mid
c\mid\mid^2=3$.
  \vskip .3 cm {\bf Proposition 1.}
  {\em The characteristic polynomial $CH_d(x)$ of the quadratic form
$2Q_d$ equals}
  $$
CH_d(x)\; =(x^3-3x+2m)(x^3-3x-2m)(x^3-3x+2n)^2
$$
\vskip .3 cm {\em Proof.} We have
  $$
  Q_d(X,Y,Z)=X^t M_c Y+X^t M_b^tZ+Y^t M_a Z
  $$
$$=\mid\mid c\mid\mid X^t O_c Y+\mid\mid b\mid\mid X^t O_b^t
Z+\mid\mid a\mid\mid Y^t O_a Z.
$$
Let us perform the orthogonal change of variables given by:
$$
x=O_c^t X,\quad y=Y,\quad
z=O_aZ.\quad\quad\quad\quad\quad\quad\quad\quad
$$
Then in these new variables the form $ Q_d$ becomes equal to
$$
\tilde Q_d(x,y,z)=\mid\mid c\mid\mid x^t     y+\mid\mid b\mid\mid
x^t O_c^t\cdot O_b^t\cdot O_a^t z+\mid\mid a\mid\mid y^t
z.\quad\quad\quad
$$
  Thus, the matrix of the form $2\tilde Q_d$ is the following block
matrix:
$$\tilde M_d=     \left(%
\begin{array}{ccc}0_4
    & \mid\mid c\mid\mid I_4 & \quad\mid\mid b\mid\mid O_{abc}^t \\
   \quad\mid\mid c\mid\mid I_4&0_4 & \mid\mid a\mid\mid I_4     \\
   \quad\quad\mid\mid b\mid\mid O_{abc} & \mid\mid a\mid\mid I_4 & 0_4
\\
\end{array}%
\right) \quad\quad\quad\quad\quad\quad\quad\quad\quad\quad$$
  where $0_4$ and $I_4$ are the zero and the unit 4x4 matrices,
respectively,
\\$O_{abc}=O_a\cdot O_b\cdot O_c $ as above.

\par Let now $\lambda\in Sp(\tilde M_d)$, $v _\lambda=
(x_\lambda,y_\lambda,z_\lambda)$ being a corresponding
eigenvector, normalized by the conditions  $ \mid\mid v _\lambda
\mid\mid= \sqrt 3$, $ (v _\lambda,d)\ge 0$. \\The condition
$\tilde M_d\cdot v_\lambda=\lambda v_\lambda$ gives
$$
\lambda x_\lambda=\mid\mid c\mid\mid y_\lambda+\mid\mid b\mid\mid
O_{abc}^t z_\lambda\quad\quad\quad\quad\quad\quad\quad\quad
$$
$$
\lambda y_\lambda=\mid\mid c\mid\mid x_\lambda\quad+\quad\mid\mid
a\mid\mid z_\lambda\quad\quad\quad\quad\quad\quad\quad\quad
$$
$$
\lambda z_\lambda=\mid\mid b\mid\mid O_{abc}  x_\lambda+\mid\mid
a\mid\mid  y_\lambda\quad\quad\quad\quad\quad\quad\quad\quad.
$$
Multiplying the second and the third equations by $\lambda$ and
inserting in thus obtained equations the first one one finds
$$
(\lambda^2-\mid\mid c\mid\mid^2)  x_\lambda=(\mid\mid
a\mid\mid\cdot\mid\mid c\mid\mid + \lambda \mid\mid b\mid\mid
O_{abc}^t)  z_\lambda  \quad\quad\quad\quad\quad\quad\quad\quad
$$
$$
( \lambda^2-\mid\mid a\mid\mid^2)  z_\lambda=(\mid\mid
a\mid\mid\cdot\mid\mid c\mid\mid + \lambda \mid\mid b\mid\mid
O_{abc})  x_\lambda  \quad\quad\quad\quad\quad\quad\quad\quad
  $$
which implies
  $$
  ( \lambda^2-\mid\mid a\mid\mid^2)(
\lambda^2-\mid\mid c\mid\mid^2) x_\lambda=(\mid\mid
a\mid\mid\cdot\mid\mid c\mid\mid + \lambda \mid\mid b\mid\mid
O_{abc}^t)(\mid\mid a\mid\mid\cdot\mid\mid c\mid\mid + \lambda
\mid\mid b\mid\mid O_{abc})  x_\lambda
$$
and, after simplifying,
$$
\lambda (\lambda ^3I_4-3\lambda I_4-mN_{abc})x_\lambda=0 ,
$$
since $\mid\mid a\mid\mid^2+\mid\mid b\mid\mid^2+\mid\mid
c\mid\mid^2=3$,  $ m= \mid\mid a\mid\mid\cdot \mid\mid b\mid\mid
\cdot \mid\mid c\mid\mid,$ $O_{abc}  O_{abc}^t=I_4$,
$N_{abc}=O_{abc}+ O_{abc}^t\; .$\\
  Hence, either $\lambda =0$ or
   $$
    (\lambda ^3 -3\lambda)\in m\cdot Sp(N_{abc})=\{-2m, 2m,-2n ,-2n \}.
$$
  This finishes the proof for $\lambda \neq 0$. If $\lambda =0$
  we get the conditions

$$
0=\mid\mid c\mid\mid y_\lambda+\mid\mid b\mid\mid O_{abc}^t
z_\lambda\quad\quad\quad\quad\quad\quad\quad\quad
$$
$$
0=\mid\mid c\mid\mid x_\lambda\quad+\quad\mid\mid a\mid\mid
z_\lambda\quad\quad\quad\quad\quad\quad\quad\quad
$$
$$
0=\mid\mid b\mid\mid O_{abc}  x_\lambda+\mid\mid a\mid\mid
y_\lambda\quad\quad\quad\quad\quad\quad\quad\quad.
$$
immediately  implying  that $m=0$ (since else these equations give
$x_\lambda=0$) and the formula holds for this case as well. \vskip
.3 cm {\bf Corollary 1.} {\em Define the angles $\alpha,\beta \in
[0,\pi]$  by $m=\cos\alpha, n=\cos\beta$. Then}
$$
Sp( \tilde M_d)  =\{ 2 \cos(\alpha/3+\pi k/3),
2\cos(\beta/3+\pi(2l+1) /3), 2\cos(\beta/3+\pi(2l+1) /3) \},
$$
  $k=0,1,\ldots 5,\quad l=0,1,2.$
  \vskip .3cm
  {\em Proof}.
  Indeed, if we put $\lambda=2\cos\gamma$, the equations
$\lambda ^3 -3\lambda+2m=0$, $\lambda ^3 -3\lambda-2m=0$ and
    $\lambda ^3 -3\lambda+2n=0$
  become respectively,
  $\cos(3\gamma)=-\cos\alpha$,
$\cos(3\gamma)= \cos\alpha$ and $\cos(3\gamma)=-\cos\beta$ which
implies the result.

\medskip
  \par Let us now order the eigenvalues in the
decreasing order:
$$
\lambda_1\ge\lambda_2\ge\lambda_3\ge\lambda_4\ge\lambda_5\ge\lambda_6\ge\lambda_7
\ge\lambda_8\ge\lambda_9\ge\lambda_{10}
\ge\lambda_{11}\ge\lambda_{12} .
$$
\vskip .3cm
{\bf Corollary 2.} \\
1). $2\ge\lambda_1\ge \lambda_2\ge\lambda_3\ge\lambda_4\ge 1$; \\
\noindent
  2). $-1\ge\lambda_9\ge\lambda_{10}
\ge\lambda_{11}\ge\lambda_{12}\ge -2$; \\
\noindent 3). $\lambda_1\ge \sqrt 3$; $ \lambda_{12}\le-\sqrt
3$;\\ \noindent
  4). {\em If $\lambda_1/ \lambda_3=2  $
  $($resp.   $\lambda_{12}/  \lambda_{10}=2)$
   then the polynomial\\
$ CH_d(x)=(x+2)^3(x-2)(x+1)^2(x-1)^6$
  $($resp.
$CH_d(x)=(x-2)^3(x+2)(x-1)^2(x+1)^6\:)$,  and $d=v_1$ $($resp.
$d=v_{12}\;)$ where $ v_i$ is the normalized eigenvector
corresponding to}  $\lambda_i$. \vskip .3cm {\em Proof}.
  All these
conclusions, except that concerning $v_1$ (resp. $v_{12}\;)$
follow from Corollary 1 along with the following elementary lemma:
  \vskip .3 cm
  {\bf
Lemma.  }  {\em Let $F_m(x)=(x^3-3x-2m)$ with $\mid m\mid \le 1$,
and let $ x_1\ge x_2\ge x_3$ be its roots. \\$1).$ If $\;0\le m
\le 1$ then $2\ge x_1\ge \sqrt 3$, $-1\ge x_3$, and each of the
conditions $2= x_1$, $-1=   x_3$ implies $ m=1$.\\$2).$ If $\;0\ge
m\ge-1$ then $-2\le x_3\le -\sqrt 3$, $ 1\le x_1$, and each of the
conditions $1= x_1$, $-2= x_3$ implies $m=-1$.\\ $3)$. If $ \;\mid
m \mid\le 0.75$ then} $ \mid x_1 \mid > 1.38$, $ \mid x_3  \mid >
1.38$ .

  This lemma follows from the monotonicity of $ \cos(x) $ on
  $[0,\pi]$ along with the inequalities $2\cos
  ({\arccos(.75)+2\pi\over 3})<-1.38$,  $2\cos
  ({\arccos(.75)\over 3})>1.94.$

To prove that $ d=v_1$ one notes that $\lambda_1/ \lambda_3=2$
implies $m(d)=n(d)=1$ which means the function   $P$ has an
absolute maximum at $d$,  its derivative in the direction $d$
equals 1 which means that $2Q_d(d)=2 $,  i.e. $ d=v_1$. The case
of $v_{12}$ is completely similar.

  \vskip .3cm
  {\bf Corollary 3.} {\em Define
$$ \delta = \sup_{d, \mid\mid d\mid\mid=\sqrt 3}\{\;\max
\{{\lambda_+^{\perp}(d)\over \lambda_3(d)},\quad
{\lambda_-^{\perp}(d)\over \lambda_{10}(d)}\;
\}\}\quad\quad\quad\quad\quad
$$
where
  $$
  \lambda_+^{\perp}(d)=2\sup_{v\perp d,\; \mid\mid v\mid\mid=\sqrt 3} Q_d(v),\quad
\lambda_-^{\perp}(d)=2\inf_{v\perp d,\; \mid\mid v\mid\mid=\sqrt
3} Q_d(v).
$$
Then} $\delta < 3/2. $ \vskip .3cm
  {\em Proof}. By Lemma, part 3 it is true for   $\mid n \mid\le 0.75 $
  since $2/1.38<3/2.$
   Let now  $ n \ge 0.75 $ (the case $ n \le - 0.75 $
   being symmetric). Suppose that ${\lambda_+^{\perp}(d)\over
   \lambda_3(d)}\ge 3/2$, and hence $\lambda_+^{\perp}(d)\ge 1.5$
   (since $\lambda_3(d)\ge 1 $). We will show that the conditions
   $\lambda_+^{\perp}(d)\ge 1.5$  and $ n \ge 0.75 $
   are incompatible. Indeed, define
   $$
T(x,y)= P(xd+y\sqrt 3
v_+^{\perp}(d))=t_3x^3+t_2x^2y+t_1xy^2+t_0y^3
   $$
   where $ v_+^{\perp}(d)$ is a norm $ \sqrt 3$ vector on which
   $Q_d(v)$ achieves the maximum.
   We get that $t_3=T(1,0)=P(d)=n \ge 0.75 $, $t_1=T_x(0,1)=3
\lambda_+^{\perp}(d)/2 \ge
   9/4$. For any $ (x,y)$ with $x^2+y^2=1 $,  $\mid T(x,y)\mid=
    \mid P(xd+y\sqrt 3 v_+^{\perp}(d))\mid \le 1. $ Let now
    $x_0={1 \over \sqrt2 }$,  $y_0=\pm {1 \over \sqrt2 } $ where the
    sign of $y_0$ is chosen from the condition that
    $y_0{t_0+t_2 \over \sqrt2 }\ge 0$. Then   $\mid T(x_0,y_0)\mid\ge
     {t_3+t_1 \over2 \sqrt2 }\ge {3 \over 2 \sqrt2 }> 1$ which is a
contradiction.
\medskip

The following result will be used in Section 4 to deduce our main
result.
     \vskip .3cm {\bf Corollary 4.} {\em Let $u\neq v\in S_{\sqrt
3}^{11} $ be two
      vectors of norm ${\sqrt 3}$. Then
     $$ P(u)-P(v) \le 3 \sqrt 3\lambda_3(d)\mid u-v\mid /4.\hskip 2 cm
     $$
     $$P(u)-P(v) \ge 3\sqrt 3\lambda_{10}(d)\mid u-v\mid /4.\hskip 2 cm $$
where}
     $ d=\sqrt 3 {u-v\over \mid u-v\mid }.$

     {\em Proof}. Denote  $ s=\mid u-v\mid /2, $
     $z=d^{\perp}\bigcap [u,v] .$  Writing the Taylor
     development for the (cubic) function $P$, we get
     $$P(u)-P(v)=2s(Q_d(z)/\sqrt 3+ s^2 P_{ddd}/18\sqrt 3).
     $$
Since
  $$ Q_d(z)\le(1-s^2/3) \sup_{x\perp d, \mid\mid x\mid\mid= \sqrt 3}
  Q_d(x)\le 3(1-s^2/3)\lambda_+^{\perp}(d)/2\le 9(1-s^2/3)\lambda_3(d)/4,
\quad
$$
$P_{ddd}\le  2/\sqrt 3\le 2\lambda_3(d) /\sqrt 3$ we get
$$
P(u)-P(v)\le 2s\lambda_3(d)(9(1-s^2/3) /4\sqrt 3+ s^2 /27)\le
3\sqrt 3\lambda_3(d)\mid u-v\mid /4.
$$ The proof of the second inequality is completely similar.

\vskip .3 cm {\em Remark.}\hskip .5 cm
  Let us resume the spectral properties of $ 2Q_d$  when $
  d$ varies over $S=S_{\sqrt 3}^{11}. $ We have a stratification $
  V_0\subset S\supset T \supset V=V_+\bigcup V_-$ where
  $T=S_{1}^{3}\times S_{1}^{3}\times S_{1}^{3}$ is defined by the
condition
  $ m(d)=1,$ $V_+ $ (resp. $V_-$ ) is defined by $P(d)=1$ (resp.
  $P(d)=-1$), $V_0=\{d: m(d)=0 \}$ ; each of $V_+ $ and $V_- $ is
diffeomorphic to
  $S_{1}^{3}\times S_{1}^{3}.$ On $V_+ $ (resp. $V_-$ )
  we have the characteristic polynomial $ (x+2)^3(x-2)(x+1)^2(x-1)^6$
  (resp. $(x-2)^3(x+2)(x-1)^2(x+1)^6$);  on $ S\setminus (T\bigcup
  V_0)$ we have $\sqrt 3 <\lambda_1(d)< 2,\:
1<\lambda_4(d),\:-1>\lambda_9(d),\:-\sqrt 3
   >\lambda_{12}(d)> -2.$ Finally, on $ V_0$ the polynomial equals
   $ x^4(x^2-3)^4.$

\section{Function w and map H}

In this section we show that  the function
$$ w(x)=P(x)/\mid x \mid $$
is what we want, i. e. the map
$$ H: S^{11}_1 \longrightarrow Q,\quad
H(a)=Hess(w(a))\quad\quad\quad\quad\quad\quad\quad\quad\quad\quad\quad\quad
$$
verifies the condition  (H) of Section 2.
  \vskip .3cm
  {\bf Proposition 2.} { Let $a\neq b\in S_1^{11} $. Then there
  exist two vectors $e,f\in S_1^{11} $ such that}
  $$w_{ee}(a)-w_{ee}(b) \ge \quad \mid a-b\mid /4\sqrt 3,$$
  $$w_{ff}(a)-w_{ff}(b) \le -\mid a-b\mid /4\sqrt 3.$$
  \vskip .3cm {\em Proof}. Let  $ d=\sqrt 3
  {a-b\over \mid a-b\mid}$. Recall that we denote by $ v_i,\; i=1,\ldots, 12$
  the normalized  eigenvectors of the form $2Q_d$ from Section 3
  (the eigenvalues $ \lambda_i=\lambda_i(d)$ being ordered in the decreasing order).
   Let $ V^+$ be the 3-dimensional space
  generated by $ v_1,v_2,v_3$ and let $e\in S_1^{11}\bigcap  V^+\bigcap
a^{\perp}
  \bigcap  b^{\perp} $. It means in particular that $2Q_d(e)\ge
\lambda_3(d)$. The conditions $ b\perp e$,  $ a\perp e$
  imply
  $$w_{ee}(a)=P_{ee}(a)-P(a), \;w_{ee}(b)=P_{ee}(b)-P(b),
  $$
  hence
  $$w_{ee}(a)-w_{ee}(b)=P_{ee}(a)- P_{ee}(b)-(P(a)-P(b)).
  $$
Since $ P_{ee}(x)$ is a linear function we get
$$P_{ee}(a)- P_{ee}(b)=\mid a-b\mid P_{eed}/\sqrt 3=2\mid a-b\mid Q_d(e)/\sqrt 3\ge
\lambda_3(d)\mid a-b\mid/\sqrt 3.
$$
Now, $$ P(a)-P(b)={(P(a{\sqrt 3})-P(b{\sqrt 3}))\over 3\sqrt 3}\le
{3\sqrt 3\lambda_3(d){\sqrt 3}\mid a-b\mid \over 12\sqrt 3}={
\sqrt 3\lambda_3(d)\mid a-b\mid \over 4}$$ by Corollary 4, and we
get
$$w_{ee}(a)- w_{ee}(b)\ge\lambda_3(d)\mid a-b\mid({1\over \sqrt 3}-{\sqrt 3\over 4})\ge
{\mid a-b\mid\over 4\sqrt 3} \; .
$$
The second inequality  is proven replacing $e $ by $f\in
S^{11}\bigcap  V^-\bigcap a^{\perp}
  \bigcap  b^{\perp}$ where $ V^-$ is
  generated by $ v_{10},v_{11},v_{12}.$

\bigskip

 {\bf Corollary 5.}
\smallskip

 {\em $1)$. The map $H: S^{11}_1
\longrightarrow Q,\quad H(a)=Hess(w(a))$ is a smooth embedding.
\smallskip

  $2)$. Let for $a\neq b\in S_1^{11} $
$$
\mu_1\ge\mu_2\ge\ldots
\ge\mu_{11}\ge\mu_{12}\quad\quad\quad\quad\quad\quad
$$
be the eigenvalues of $ Hess(w(a))-Hess(w(b)).$ Then}
$$ M^{-1}={1\over 1536\sqrt 3} \le - {\mu_1\over\mu_{12}} \le {1536\sqrt 3 }=M.
\quad\quad\quad\quad\quad\quad\quad
$$
\vskip .3cm {\em Proof}. \\
1). This follows immediately from Proposition 2. \\
2). An easy calculation  shows that $\mid w_{efg}(x)\mid\le 32$
for any  $e,f,g, x  \in S_1^{11}$. Hence
$$ \mid w_{ef}(a)- w_{ef}(b)\mid\le \mid w_{efd'}(d')\mid\cdot \mid
a-b\mid\le 32\mid a-b\mid
$$ for $d'=d/\sqrt 3 $.
  \noindent Since all elements  of the matrix $ Hess(w(a))-Hess(w(b))$
are of absolute value $\le 32\mid a-b\mid$, all its eigenvalues
are of absolute value $\le 12\cdot 32\mid a-b\mid.$ Using the
inequalities of Proposition 2  we get the conclusion.
\bigskip\bigskip

\centerline{REFERENCES}

\bigskip
  [A] A.D. Alexandrroff; {\it Sur les th\'eor\`emes d'unicite
pour les surfaces ferm\'ees}, Dokl. Acad. Nauk 22 (1939), 99-102.

\medskip
[CC] L. Caffarelli, X. Cabr\'e; {\it Fully Nonlinear Elliptic
Equations}, Amer. Math. Soc., Providence, R.I., 1995.

\medskip
  [CIL]  M.G. Crandall, H. Ishii, P-L. Lions; {\it User's guide
to viscosity solutions of second order partial differential
equations}, Bull. Amer. Math. Soc. (N.S.)  27(1) (1992), 1-67.

\medskip
[GT] D. Gilbarg, N. Trudinger; {\it Elliptic Partial Differential
Equations of Second Order, 2nd ed.}, Springer-Verlag,
Berlin-Heidelberg-New York-Tokyo, 1983.

\medskip
[HNY] Q. Han, N. Nadirashvili, Y. Yuan, {\it Linearity of
homogeneous order-one solutions to elliptic equations in dimension
three}, Comm. Pure Appl. Math. 56 (2003), 425-432

\medskip

[K] N.V. Krylov; {\it Nonlinear Elliptic and Parabolic Equations
of Second Order}, Reidel, 1987.
\medskip

\medskip
[NY] N. Nadirashvili, Y. Yuan, {\it Homogeneous solutions to fully
nonlinear elliptic equation}, Proc. AMS  134 (2006), no 6.

\medskip
[N] L. Nirenberg; {\it On Nonlinear Elliptic Partial Differential
Equations and H\"older Continuity}, Comm. Pure Appl. Math. 6
(1953), 103--156.

\end{document}